\documentclass[12pt,leqno]{amsart}
\usepackage{amssymb}
\usepackage{amsmath,amssymb,color}
\theoremstyle{plain}% Theorem-like structures provided by amsthm.sty
\newtheorem{theorem}{Theorem}[section]
\newtheorem{lemma}[theorem]{Lemma}
\newtheorem{corollary}[theorem]{Corollary}
\newtheorem{proposition}[theorem]{Proposition}

\theoremstyle{definition}

\newtheorem{example}[theorem]{Example}

\theoremstyle{remark}

\numberwithin{equation}{section}
\newcommand{\sump}{\sum_{i=1}^n}
\newcommand{\sumq}{\sum_{j=1}^m}
\newcommand{\qqqt}{\ooo{Q}_T}
\newcommand{\ssst}{\ooo{S}_T}

\newcommand{\ep}{\varepsilon}

\newcommand{\va}{\varphi}
\newcommand{\ppp}{\partial}

\newcommand{\pppa}{\partial_t^{\alpha}}

\newcommand{\R}{\mathbb{R}}

\newcommand{\ooo}{\overline}
\newcommand{\OOO}{\Omega}

\newcommand{\hhalf}{\frac{1}{2}}

\allowdisplaybreaks
\title
[]
{
On the maximum principle for the multi-term fractional transport equation
}

\author{Yuri Luchko$^1$,  Anna Suzuki$^2$,
Masahiro Yamamoto$^{3,4,5}$}

\thanks{
$^1$ Department of Mathematics, Physics, and Chemistry, Beuth Technical University of Applied Sciences Berlin, Luxemburger Str. 10, 13353 Berlin, Germany; e-mail: {\tt luchko@beuth-hochschule.de} \\
$^2$ Institute of Fluid Science, Tohoku University, 2-1-1 Katahira, Aoba-ku, Sendai, Miyagi 980-8577, Japan; e-mail: {\tt anna.suzuki@tohoku.ac.jp} \\
$^3$ Graduate School of Mathematical Sciences, The University
of Tokyo, Komaba, Meguro, Tokyo 153-8914, Japan \\
$^4$ Honorary Member of Academy of Romanian Scientists, 
Splaiul Independentei Street, no 54,
050094 Bucharest Romania \\
$^5$ Peoples' Friendship University of Russia 
(RUDN University) 6 Miklukho-Maklaya St, Moscow, 117198, Russian Federation;
e-mail: {\tt myama@ms.u-tokyo.ac.jp}
}

\date{}
\begin{document}
\maketitle

\baselineskip 18pt

\begin{abstract}
In this paper, we prove a maximum principle for the general multi-term space-time-fractional transport equation and 
apply it for establishing  uniqueness of solution to an 
initial-boundary-value problem for this equation. We also derive some comparison principles for solutions to the initial-boundary-value problems with different problem data.
Finally, we present a maximum principle for the Cauchy problem for a time-fractional transport equation on an unbounded domain. 
%35R30, 35R11
\end{abstract}

\section{Introduction}
\label{s1}

Within the last few decades, fractional calculus in general and fractional partial differential equations became a very popular and important topic both in mathematics and in numerous applications. The framework of fractional calculus has been widely employed to describe several physical phenomena including anomalous diffusion and anomalous transport processes in various areas, such as material science \cite{MN}, medical engineering \cite{LX,W}, electrical engineering \cite{SG}, hydrology \cite{BMR}, geological engineering \cite{OAHA,SNFCH4}, and the earth systems \cite{ZSSZH}.  

One of the most investigated and used fractional partial differential equations is the time-fractional
diffusion equation. In the one-dimensional case and on the finite space- and time intervals, the time-fractional
diffusion equation with the convection and reaction terms is formulated as follows:
$$
\pppa u(x,t) = \ppp_x^2u(x,t) - q(x,t)\ppp_xu(x,t) + r(x,t)u(x,t),\ 0<\alpha \le 1,\ 
\ 0<x<\ell, \ 0<t<T.                 \eqno{(1.1)}
$$
% In this paper, we mainly address the one-dimensional case and the equations defined on the finite space- and time intervals.
For $0<\alpha<1$, by $\pppa$ we denote the Caputo fractional derivative  (\cite{Po}):
$$
\pppa u(x,t) = \frac{1}{\Gamma(1-\alpha)}\int^t_0
(t-s)^{-\alpha}\, \frac{\partial}{\partial s} u(x,s)\, ds,
$$
where $\frac{\partial}{\partial s} u(\cdot,s)$ is assumed to belong to the space $L^1(0,T)$ for any $x\in (0,\ell)$ and
$\Gamma$ is the Euler gamma function. For $\alpha =1$, $\pppa$ is interpreted as the conventional first order derivative. 
As usual, we set $\ppp_x := \frac{\ppp}{\ppp x}$ and 
$\ppp_x^2 := \frac{\ppp^2}{\ppp x^2}$.

The initial-boundary-value problems for the equation (1.1) in different settings and properties of their solutions
have been already intensively studied in the literature.  Especially for the unique existence 
of solutions to the initial-boundary-value problems for the equation (1.1) and its generalizations, we refer to
\cite{KiYa,KRY,Lu10,Lu2,SY,Za2} to mention only few of many relevant publications. Moreover, there are many
works on numerical analysis of the fractional 
partial differential equations, but our focus in this paper is on their analytical treatment and we do not refer to any publications regarding numerical methods.

In this paper, we address the following time-fractional 
%(and/or spatial-fractional) 
transport equation with the Caputo fractional derivative of the order $\alpha,\ 0<\alpha < 1$:
$$
p(x,t)\pppa u(x,t) + q(x,t)\ppp_xu(x,t) = r(x,t)u(x,t) + F(x,t),
\quad 0<x<\ell, \, 0<t<T                 \eqno{(1.2)}
$$
along with the boundary and initial conditions
$$
u(0,t) = g(t), \quad 0<t<T,               \eqno{(1.3)}
$$
$$
u(x,0) = a(x), \quad 0<x<\ell,      \eqno{(1.4)}
$$
respectively, as well as its multi-term time-space-fractional generalizations
which we formulate in Section 3. 

In what follows, we assume  the inclusions 
$$
p, q, r\in C([0,\ell]\times [0,T]) \quad
                                               \eqno{(1.5)}
$$
as well as some conditions on the signs of the functions $p, q, r$ that 
we formulate in due time. 
Taking into 
consideration the outgoing and ingoing sub-boundaries in the case
of the transport equation (1.2) with $\alpha=1$, it is natural to prescribe the boundary condition 
(1.3) at the point $x=0$, not at the point $x=\ell$.

Throughout the paper, we assume the 
existence of a solution to the initial-boundary-value problem 
(1.2)-(1.4) that satisfies the following inclusions:
$$
u\in C([0,\ell] \times [0,T]), \quad u(\cdot,t) \in W^{1,1}(0,\ell),
\quad u(x,\cdot) \in W^{1,1}(0,T),        \eqno{(1.6)}
$$
where $W^{1,1}(0,T) = \{ g;\, g, \, \ppp_tg \in L^1(0,T)\}$ and 
$W^{1,1}(0,\ell) = \{ a;\, a, \, \ppp_xa \in L^1(0,\ell)\}$.
 
The time-fractional transport equation (1.2) was already employed  
for modeling various anomalous transport processes including the mass and heat transfer for characterizing geothermal reservoirs 
(\cite{SNFCH,SNFCH2,SNFCH3,SHLH}).
However, compared to the comprehensive results already obtained for the time-fractional diffusion equation of type (1.1), it may be 
a surprise that until now only few theoretical publications were devoted to the fractional transport equations, for instance, to  the problem of
 unique existence of solution to the initial-boundary-value problem
(1.2)-(1.4).  For a treatment of the viscosity solutions to the time-fractional 
transport equations we refer 
to \cite{N}.

For $\alpha=1$, the equation (1.2) is the classical and well studied
transport equation. It is well known that the solutions to 
the initial-boundary-value problems of type (1.2)-(1.4) with $\alpha=1$ can be analyzed by the method of characteristics.
However, for $0< \alpha < 1$, the method of characteristics does not work
and thus the properties of the fractional transport equation (1.2) are not yet well investigated.  

Another important approach to analysis of the solution properties to the partial differential equations is the maximum principle (\cite{Pro}).   
For the multi-dimensional time-fractional diffusion equation of type (1.1), the maximum 
principle in different settings has been proved in \cite{Lu1,LuY1,Za1} for the case of
the Caputo time-fractional derivative,   
in  \cite{ARL} for the case of the Riemann-Liouville derivative, and in  \cite{LuY0} for the case of the general fractional derivative introduced in  \cite{Koch11}. For more results regarding the maximum principles for the fractional partial differential equations we refer to the 
surveys  \cite{LuY2,LuY3}.

However, to the best knowledge of the authors, no maximum principle for the fractional transport equations has 
 been yet established. In this paper, we formulate and prove a maximum principle for the multi-term space-time-fractional transport equation and derive some of its useful consequences.   Since the method of characteristics does not work
for the fractional transport equations,   the maximum and comparison principles are worth employing as an alternative methodology for their analytical treatment.

The rest of the paper is organized as follows. In Section \ref{s1-2}, we start with a simple case of the fractional transport equation (1.2) and illustrate the main ideas behind the derivations in the general case. In Section \ref{s2}, our main results are formulated. The next two sections are devoted to the proofs of two main theorems stated in   Section \ref{s2}. Finally, in the last section, we provide some concluding remarks and directions for further research. 

\section{Illustrating example}
\label{s1-2}

Before stating and proving our main results, in this section, we address the following simple particular case of the time-fractional transport equation (1.2):
$$
\pppa u(x,t)  + q_0\ppp_xu(x,t) = r(x,t)u(x,t), \ 0<\alpha \le 1,\ 
\quad 0<x<\ell, \, 0<t<T               \eqno{(1.7)}
$$
along with the boundary condition (1.3) and the initial condition (1.4).  In the equation (1.7),  $q_0>0$ is a constant and we assume that the condition
$$
r(x,t) < 0, \quad 0\le x \le \ell, \, 0\le t \le T
                                               \eqno{(1.8)}
$$
is satisfied. It is worth mentioning that the condition (1.8) can be replaced with a weaker condition $r(x,t)\le 0$ on $[0,\ell]\times [0,T]$. However, for simplicity of the proofs, in this section we suppose that the stronger inequality (1.8) is satisfied. 
For the solution $u=u(x,t)$ of the initial-boundary-value problem (1.7), (1.3) and (1.4), the following result holds true:

%{\bf Proposition 1.}\\
%{\it
\begin{proposition}
\label{p1}
Let $u=u(x,t)$ satisfy the inclusion (1.6) and 
$u(0,t) \ge 0$ for $0\le t \le T$ and $u(x,0) \ge 0$ for 
$0\le x \le \ell$.  Then $u=u(x,t)$ is non-negative on the whole domain $\ooo{Q}_T:= [0,\ell]\times [0,T]$, i.e., 
$$
u(x,t) \ge 0, \quad  0\le x \le \ell, \, 0\le t \le T.
$$
%}
%\\
\end{proposition}

\begin{proof}
The proof of the proposition essentially relies on the extremum principle for the Caputo
fractional derivative.
%\\
%{\bf Lemma 1 (Luchko \cite{Lu1}).}\\
%{\it
\begin{lemma}[\cite{Lu1}]
\label{l1}
Let $f \in C[0,T] \cap W^{1,1}(0,T)$ attain
its maximum (resp. its minimum) over the closed interval $[0,T]$ at a point $t_0 
\in (0,T]$.  Then for any $\alpha\in (0,1]$ the inequality
$$
(\pppa f)(t_0) \ge 0 \quad \mbox{(resp. $(\pppa f)(t_0) \le 0$)}
$$
holds true.
\end{lemma}
%}
%\\
%Proof of Proposition 1.\\
We employ an indirect proof and assume that the conclusion of the proposition does not hold true.  Then there exists a point 
$(x_0,t_0) \in \Omega_T$ such that 
$$
u(x_0,t_0) := \min_{(x,t)\in \Omega_T} u(x,t) < 0.
$$
Since $u(0,t) \ge 0$ and $u(x,0) \ge 0$, we conclude that 
$x_0>0$ and  $t_0 > 0$.  By Lemma \ref{l1}, we have $\pppa u(x_0,t_0) \le 0$.
Moreover, since $u(x,t_0)$ as a function in $x$ attains 
its minimum  at the point $x=x_0>0$, we get the inequality $\ppp_xu(x_0,t_0) \le 0$.  Hence
$$
0\ge \pppa u(x_0,t_0) = -q_0\ppp_xu(x_0,t_0) + r(x_0,t_0)u(x_0,t_0)
\ge r(x_0,t_0)u(x_0,t_0) > 0
$$ 
because $r(x_0,t_0) < 0$ and $u(x_0,t_0) < 0$.
With the last inequality, we arrived to a contradiction and the proof of Proposition \ref{p1} is completed.
\end{proof}
%\vspace{0.3cm}

Proposition \ref{p1} is quite preliminary and serves just for illustration of our method. In 
the next section, we present a maximum principle for the more general multi-term 
time-space-fractional 
 transport equation that is valid under the weaker conditions on the problem data compared to the ones formulated  in Proposition \ref{p1}.
\section{Main results}
\label{s2}

In this section, we  address an initial-boundary-value problem for a one-dimensional multi-term  time-space-fractional transport equation defined on the bounded domain  
$
Q_T := (0,\ell) \times (0,T),\ \ell >0,\ T>0$ with the boundary  
$S_T := \{ (x,0);\, 0<x<\ell\} \cup \{ (0,t);\, 0<t<T\}$.

To formulate the equation, we first introduce the functions $p_i=p_i(x,t),\ i=1,\dots,n$ and $q_j=q_j(x,t),\ j=1,\dots,m $ and the constants $\alpha_i, \ i=1,\dots,n$ and  $\beta_j,\ j=1,\dots,m $ that satisfy the following conditions and inclusions: 
$$
\left\{ \begin{array}{rl}
&0 < \alpha_1 < \cdots < \alpha_n \le 1 ,\quad
0 < \beta_1 < \cdots < \beta_m \le 1, \\
& p_i \in C(\ooo{Q}_T), \, q_j\in C(\ooo{Q}_T) \quad \mbox{ for }
1\le i \le n, 1\le j \le m, \\
& p_i(x,t) \ge 0, \, 1\le i \le n, \quad
q_j(x,t) \ge 0, \, 1\le j \le m, \quad (x,t) \in Q_T,\\
& \sum_{i=1}^{n} p_i(x,t) > 0 \quad \mbox{for all } (x,t) \in 
\ooo{Q}_T.
\end{array}\right.
                                           \eqno{(2.1)}
$$

For the given functions $r, F \in C(\ooo{Q}_T)$, the one-dimensional multi-term  time-space-fractional transport equation is introduced as follows:
$$
\sump p_i(x,t)\ppp_t^{\alpha_i}u(x,t)
+ \sumq q_j(x,t)\ppp_x^{\beta_j}u(x,t)
= r(x,t)u(x,t) + F(x,t),  (x,t)\in Q_T,  \eqno{(2.2)}
$$
where for $0<\beta < 1$ the space-fractional Caputo derivative is defined by the formula 
$$
\ppp_x^{\beta}u(x,t) = \frac{1}{\Gamma(1-\beta)}
\int^x_0 (x-y)^{-\beta} \frac{\partial}{\partial y}u(y,t) dy 
$$
in analogy to the time-fractional Caputo derivative. Of course, for $\beta=1$ the Caputo fractional derivative is interpreted as the conventional first order derivative.  

In what follows, we always assume that the  function $r=r(x,t)$ is non-positive, i.e.,
$$
r(x,t) \le 0, \quad (x,t) \in \ooo{Q_T}     \eqno{(2.3)}
$$
and that any solution $u=u(x,t)$ to the equation (2.2) satisfies the regularity conditions (1.6).

Our main results are formulated in the following two theorems.

\begin{theorem}
\label{t1}
% Theorem 1.}\\
%{\it
(i) Let $F(x,t)\le 0$ for $(x,t)\in\ooo{Q}_T$.  Then
$$
\max_{(x,t)\in \ooo{Q}_T} u(x,t) \le \max \{0, \, \max_{(x,t)\in \ssst} u(x,t)\}.  \eqno{(2.4)}
$$
\\
(ii) Let $F(x,t)\ge 0$ for $(x,t)\in \ooo{Q}_T$.  Then
$$
\min_{(x,t)\in \ooo{Q}_T} u(x,t) \ge \min \{0, \, \min_{(x,t)\in \ssst} u(x,t)\}.  \eqno{(2.5)}
$$
In the case of $r\equiv 0$ in $Q_T$, the inequalities (2.4) and (2.5) can be replaced by the equalities 
$$
\max_{(x,t)\in \ooo{Q}_T} u(x,t) = \max_{(x,t)\in \ssst} u(x,t)       \eqno{(2.4)'}
$$
and 
$$
\min_{(x,t)\in \ooo{Q}_T} u(x,t) = \min_{(x,t)\in \ssst} u(x,t),  \eqno{(2.5)'}
$$
respectively.
%}
%\\
\end{theorem}

From Theorem \ref{t1}, we readily derive the following useful consequence:
%{\bf Corollary 1.}\\
%{\it
\begin{corollary}
\label{c1}
If $F\equiv 0$ in $Q_T$, then 
$$
\max_{(x,t)\in \ooo{Q}_T} u(x,t) = \max_{(x,t)\in \ssst} u(x,t), \quad 
\min_{(x,t)\in \ooo{Q}_T} u(x,t) = \min_{(x,t)\in \ssst} u(x,t).         \eqno{(2.6)}
$$
\end{corollary}

In its turn, this corollary immediately yields an uniqueness result.
\begin{corollary}[uniqueness of solution]
\label{c2}
Let the functions $u_1=u_1(x,t)$ and $u_2=u_2(x,t)$ satisfy the equation (2.2) and the regularity conditions (1.6).  If $u_1(x,t) = u_2(x,t)$ 
on the boundary $S_T$ of the domain $Q_T$, then $u_1(x,t) = u_2(x,t)$ on the whole domain $Q_T$.
\end{corollary}
\begin{proof}
Indeed, setting $u:= u_1 - u_2$, we see that the function
$u$ satisfies (1.6) and (2.2) with $F\equiv 0$.  Therefore Corollary \ref{c1}
implies $\max_{(x,t)\in \qqqt} u(x,t) = \min_{(x,t)\in \ssst} u(x,t) = 0$, which  means that 
$u_1 \equiv  u_2$ in $Q_T$.
\end{proof}

Moreover, Theorem \ref{t1} implicates some important comparison principles.
Let the function $u=u_{r,a,g,F}(x,t)$ satisfy the regularity conditions (1.6), the equation (2.2), and the following initial and boundary conditions:
$$
u(0,t) = g(t), \quad 0<t<T, \quad u(x,0) = a(x), \quad 0<x<\ell.
                                     \eqno{(2.7)}
$$
\begin{corollary}[comparison principles]
\label{c3}

(i) Let $F_1(x,t) \ge F_2(x,t)$ for $(x,t)\in \qqqt$, $g_1(t) \ge g_2(t)$ for $t\in [0,T]$,
and $a_1(x) \ge a_2(x)$ for $x\in [0, \ell]$.  Then
$$
u_{r,a_1,g_1,F_1}(x,t) \ge u_{r,a_2,g_2,F_2}(x,t), \  (x,t) \in \qqqt.
$$

(ii)  Let $F_1(x,t) \ge F_2(x,t)\ge 0$ for $(x,t)\in \qqqt$, $g_1(t) \ge g_2(t)\ge 0$ for $t\in [0,T]$,
and $a_1(x) \ge a_2(x)\ge 0$ for $x\in [0, \ell]$.   If $0\ge r_1(x,t) \ge r_2(x,t)$ for
$(x,t) \in \qqqt$, then
$$
u_{r_1,a_1,g_1,F_1}(x,t) \ge u_{r_1,a_2,g_2,F_2}(x,t), \  (x,t) \in \qqqt.
$$
\end{corollary}

For the similar comparison principles for the time-fractional
diffusion equation we refer the readers to \cite{LuY1}.

Now we formulate a maximum principle for the following Cauchy problem for a time-fractional transport equation of order $\alpha$, $0<\alpha<1$ defined on an unbounded domain $\Omega_T := \R \times (0,T)$:
$$
\pppa u(x,t) + q(x)\ppp_xu(x,t) = F(x,t), \quad (x,t) \in \Omega_T 
                                                             \eqno{(2.8)}
$$
along with the initial condition
$$
u(x,0) = a(x), \quad x\in \R.                       \eqno{(2.9)}
$$
In (2.9), we assume the inclusion $a \in W^{1,1}_{loc}(\R)$ that means that $a\vert_{(-X,X)} \in W^{1,1}(-X,X)$ for 
any $X>0$. Evidently, equation (2.8) is a particular case of the multi-term  time-space-fractional transport equation (2.2). The following maximum principle is valid:

\begin{theorem}
\label{t2}
Let $u \in C(\R \times [0,T]) \cap L^{\infty}(\R \times (0,T))$ satisfy 
(2.8) and (2.9) and the inclusions
$$
u(x,\cdot) \in W^{1,1}(0,T), \quad u(\cdot,t) \in W^{1,1}_{loc}(\R).
$$
Moreover, we assume
$$
\int^0_{-\infty} \left\vert \frac{1}{q(\xi)}\right\vert d\xi 
= \infty.                 \eqno{(2.10)}
$$
If  $F(x,t)\le 0$ for $(x,t)\in Q_T$, then
$$
\sup_{(x,t) \in \Omega_T } u(x,t) = \sup_{x\in \R} a(x).                   \eqno{(2.11)}
$$
If  $F(x,t)\ge 0$ for $(x,t)\in Q_T$, then
$$
\inf_{(x,t) \in \Omega_T } u(x,t) = \inf_{x\in \R} a(x).       \eqno{(2.12)}
$$
\end{theorem}

To demonstrate the statement of Theorem \ref{t2}, we consider two  
simple examples and address the case when the solution $u$ can be represented in the form $u(x,t) = \psi(t) + \va(x)$ for 
$x\in \R$ and $0<t<T$, where $\psi \in C^1[0,T]$ and 
$\va \in C^1(\R) \cap L^{\infty}(\R)$.

\begin{example}
Let us suppose that $\frac{d\psi}{dt}(t) \le 0$ for $0<t<T$.  Then we can easily verify
that $\pppa \psi(t) \le 0$ for $0 < t < T$.  If we choose $\va(x)$ such that
$\frac{d\va}{dx}(x) \le 0$ for $x\in \R$, then $u(x,t) = \psi(t) 
+ \va(x)$ satisfies (2.8) with $q(x,t)\ge 0$ and $F(x,t)= \pppa \psi + q(x)\ppp_x \va \le 0$.  
Hence the equality (2.11) holds true.  However, in this case, (2.11) is trivial because
$\frac{\ppp u}{\ppp t}(x,t) = \frac{d\psi}{dt}(t) \le 0$ for 
$0\le t\le T$ and so we immediately see that $u(x,t) \le u(x,0)$
for $x\in \R$ and $0<t<T$.
\end{example}
\begin{example}
Theorem \ref{t2} is less trivial if the inequalities $\frac{d\psi}{dt}(t) \le 0$ and 
$\frac{d\va}{dx}(x) \le 0$ do not hold.  
For example, let $\va(x) = e^{-x^2}$, $x\in \R$.  
Then $\frac{d\va}{dx}(x) > 0$ for $x<0$.
 
If $\pppa \psi(t) \le -\sqrt{2}e^{-\hhalf}
= -\max_{x\in\R} \left\vert \frac{d\va}{dx}(x) \right\vert$, then
$$
\pppa u(x,t) + \ppp_xu(x,t)
= \pppa \psi(t) + \frac{d\va}{dx}(x) 
\le -\max_{x\in\R} \left\vert \frac{d\va}{dx}(x) \right\vert
+ \frac{d\va}{dx}(x) \le 0, \quad x\in \R, \, 0\le t \le T,
$$
i.e., the equation (2.8) holds valid with a function $F=F(x,t)\le 0$ and $q(x) = 1$, $x \in \R$.
The statement of Theorem \ref{t2} is that  the inequality 
$\pppa \psi(t) \le -\sqrt{2}e^{-\hhalf}$ for $0<t<T$ implies the inequality 
$\frac{d\psi}{dt}(t) \le 0$ for $0<t<T$ that is not trivial.
\end{example}
 
\section{Proof of Theorem \ref{t1}}
\label{s3}

In this section, we present a proof of Theorem \ref{t1} that is based on Lemma \ref{l1} and carried out similarly to the proof of Theorem 2  
from \cite{Lu1}.  It suffices to prove the inequality (2.4) and the equality (2.4)' because the inequality (2.5) and the equality (2.5)' can be proved
by replacing $u$ by $-u$ and arguing in the same way.

\begin{proof}
We prove the inequality (2.4) by contradiction.
Assume that (2.4) does not hold true.
Then there exist $x_0 \in [0,\ell]$ and $t_0 \in [0,T]$ such that 
$$
u(x_0,t_0) > M:= \max\{ 0, \max_{(x,t)\in \ssst} u(x,t)\} \ge 0.
$$
Now we set
$$
\ep := u(x_0,t_0) - M > 0                 \eqno{(3.1)}
$$
and introduce an auxiliary function $w(x,t)$, which is the same as the one employed 
in \cite{Lu1}:
$$
w(x,t):= u(x,t) + \frac{\ep}{2T}(T-t), \ (x,t) \in Q_T.
%                                                      \eqno{(3.2)}
$$
It is easy to calculate that 
$$
\ppp_t^{\alpha_i}w(x,t) = \ppp_t^{\alpha_i}u(x,t)
- \frac{\ep t^{1-\alpha_i}}{2T\Gamma(2-\alpha_i)}, \quad 
1\le i \le n.
$$
Therefore we have the following equality
$$
\sump p_i(x,t)\ppp_t^{\alpha_i}w(x,t) 
+ \sumq q_j(x,t)\ppp_x^{\beta_j}w(x,t)
$$
$$
= r(x,t)\left( w(x,t) - \frac{\ep}{2T}(T-t)\right)
- \frac{\ep}{2T}\sump 
\frac{p_i(x,t)t^{1-\alpha_i}}{\Gamma(2-\alpha_i)} + F(x,t),\ (x,t) \in  Q_T.                     \eqno{(3.2)}
$$
By definition of $w$, the inequality
$$
w(x,t) \ge u(x,t), \quad (x,t) \in \qqqt
$$
holds true. On the other hand, the condition (3.1) yields
$$
w(x_0,t_0) \ge u(x_0,t_0) = M + \ep.     \eqno{(3.3)}
$$
Since $u(x,t) \le M$ for $(x,t) \in \qqqt$, 
the chain of inequalities
\begin{align*}
& w(x_0,t_0) \ge M+\ep \ge \ep + u(x,t) 
\ge \ep + w(x,t) - \frac{\ep}{2T}(T-t)\\
\ge& \ep + w(x,t) - \frac{\ep}{2} 
= w(x,t) + \frac{\ep}{2} > w(x,t)
\end{align*}
holds true for any $(x,t) \in \ssst$ that in its turn implies
the inequality  $\max_{\ssst} w < w(x_0,t_0)$.

This means that if $w$ attains its maximum over $\qqqt$ at the point 
$(x_1, t_1)$, then
$$
(x_1, t_1) \not\in \ssst
$$
and therefore
$$
x_1 > 0, \quad t_1 > 0.      \eqno{(3.4)}
$$
Moreover, by (3.3) and $M\ge 0$, we obtain the estimates
$$
w(x_1,t_1) \ge w(x_0,t_0) \ge M+\ep \ge \ep.  \eqno{(3.5)}
$$
Because of the conditions (3.4), we may apply Lemma \ref{l1} and get the following inequalities (in the case $\alpha_i=1$ 
or $\beta_j = 1$, these inequalities are well known in calculus):
$$
\ppp_t^{\alpha_i}w(x_1,t_1) \ge 0, \quad 1\le i\le n,
\quad \ppp_x^{\beta_j}w(x_1,t_1) \ge 0, \quad 1\le j\le m.
$$
Hence
$$
\sump p_i(x_1,t_1)\ppp_t^{\alpha_i}w(x_1,t_1) 
+ \sumq q_j(x_1,t_1)\ppp_x^{\beta_j}w(x_1,t_1) \ge 0.
                                          \eqno{(3.6)}
$$
It follows from the inequality (3.5) that 
$$
w(x_1,t_1) - \frac{\ep}{2T}(T-t_1)
\ge \ep - \frac{\ep}{2T}T = \frac{\ep}{2} > 0
$$
and
$$
r(x_1,t_1)\left( w(x_1,t_1) - \frac{\ep}{2T}(T-t_1)\right)
\le 0                     \eqno{(3.7)}
$$
because of the assumption (2.3). 

Moreover, the inequality
$$
\sump \frac{p_i(x_1,t_1)t_1^{1-\alpha_i}}{\Gamma(2-\alpha_i)} > 0
$$
holds true. 
Indeed, let us assume that $\sump \frac{p_i(x_1,t_1)t_1^{1-\alpha_i}}{\Gamma(2-\alpha_i)} 
= 0$.  Since $\frac{p_i(x_1,t_1)t_1^{1-\alpha_i}}{\Gamma(2-\alpha_i)} 
\ge 0$ for $1\le i \le n$ by the assumptions (2.1), we get 
$\frac{p_i(x_1,t_1)t_1^{1-\alpha_i}}{\Gamma(2-\alpha_i)} = 0$ for any $i=1,\dots,n$, 
which implies $p_i(x_1,t_1) = 0,\ i=1,\dots,n$.
Therefore $\sump p_i(x_1,t_1) = 0$ that contradicts the last of the conditions (2.1).

Thus, the inequality 
$$
-\frac{\ep}{2T}\sump \frac{p_i(x_1,t_1)t_1^{1-\alpha_i}}
{\Gamma(2-\alpha_i)} < 0               \eqno{(3.8)}
$$
holds true.
Using the condition $F(x,t)\ge 0$ and substituting the inequalities (3.6) - (3.8) into the formula (3.2), we arrive at a 
contradiction that proves the inequality (2.4) (and hence the inequality (2.5)). 

Now we proceed with a proof of the equality (2.4)' and assume that
 $r\equiv 0$ in $Q_T$.  Then, instead of $M$ as in the previous proof, we set 
$M_0:= \max_{(x,t)\in \ssst} u(x,t)$. We repeat the same arguments as above to obtain the inequalities $x_1>0$ and
$t_1 > 0$, where $w(x_1,t_1)$ is the maximum of $w = w(x,t)$ over $\qqqt$, and the equation
\begin{align*}
& \sump p_i(x_1,t_1)\ppp_t^{\alpha_i}w(x_1,t_1)
+ \sumq q_j(x_1,t_1)\ppp_x^{\beta_j}w(x_1,t_1)\\
=& -\frac{\ep}{2T} \sump 
\frac{p_i(x_1,t_1)t_1^{1-\alpha_i}}{\Gamma(2-\alpha_i)} + F(x_1,t_1)
\end{align*}
in place of the equation (3.2).
Furthermore, we can verify that the inequalities (3.6) and (3.8)  hold true and then arrive at  a contradiction
similar to the one formulated above.  The only difference to the previous proof is that we cannot use the inequality (3.5) 
because the case $M_0 < 0$ may occur. However, we do not need it this time because  of the assumption 
$r\equiv 0$. The proof of Theorem \ref{t1} is completed.
\end{proof}

Now we prove Corollary \ref{c3}.

\begin{proof}
First we prove the part (i) of Corollary \ref{c3}.
We start by setting $u := u_{r,a_1,g_1,F_1} - u_{r,a_2,g_2,F_2}$ and 
$F:= F_1-F_2$. Then $F(x,t)\ge 0, \ (x,t)\in Q_T$ and the function $u$ satisfies the equation
$$
\sump p_i(x,t)\ppp_t^{\alpha_i}u(x,t)
+ \sumq q_j(x,t)\ppp_x^{\beta_j}u(x,t) = r(x,t)u(x,t) + F(x,t),\ (x,t)\in Q_T
$$
and the inequalities
$$
u(0,t) \ge 0, \quad 0\le t \le T, \quad u(x,0) \ge 0, \quad 
0\le x \le \ell.
$$
Thus $\min_{(x,t)\in \ssst} u(x,t) \ge 0$  and  
$\min\{ 0,\, \min_{(x,t)\in \ssst} u(x,t)\} = 0$.  Since $F(x,t)\ge 0, \ (x,t)\in Q_T$, we can apply 
the inequality (2.5) from Theorem \ref{t1} and get the inequality
$$
\min_{(x,t)\in \qqqt} u(x,t) \ge \min\{0, \, \min_{(x,t)\in\ssst} u(x,t)\} = 0.
$$
The proof of (i) is completed.

Then we proceed with a proof of the part (ii) of Corollary \ref{c3}.

Because $a_2(x) \ge 0,\ x\in[0,\ell]$ and $g_2(t)\ge 0,\ t\in [0,T]$, Theorem \ref{t1} yields
the inequality
$$
u_{r_2,a_2,g_2,F_2}(x,t) \ge 0,\ (x,t) \in \qqqt.    \eqno{(3.9)}
$$
Now we again use the notations $u := u_{r,a_1,g_1,F_1} - u_{r,a_2,g_2,F_2}$ and 
$F:= F_1-F_2$. Then $F(x,t)\ge 0, \ (x,t)\in Q_T$ and the function $u$ satisfies the equation
\begin{align*}
& \sump p_i(x,t)\ppp_t^{\alpha_i}u(x,t)
+ \sumq q_j(x,t)\ppp_x^{\beta_j}u(x,t) \\
=& r_1(x,t)u(x,t) + (r_1(x,t)-r_2(x,t))u_{r_2,a_2,g_2,F_2}(x,t) + F(x,t), \ (x,t) \in Q_T
\end{align*}
and the inequalities
$$
u(0,t) \ge 0, \quad 0\le t \le T, \quad u(x,0) \ge 0, \quad 
0\le x \le \ell.
$$
Because $r_1(x,t)-r_2(x,t) \ge 0,\ (x,t) \in Q_T$ and using the inequality (3.9), we get the inequality 
$(r_1(x,t)-r_2(x,t))u_{r_2,a_2,g_2,F_2}(x,t) + F(x,t) \ge 0,\ (x,t) \in Q_T$. Thus we can apply the inequality (2.5) to the equation for $u$ that completes the 
proof of (ii) and thus the proof of Corollary \ref{c3}.
\end{proof}

\section{Proof of Theorem \ref{t2}}
\label{s4}

\begin{proof}
The main element of our proof is a suitably chosen 
auxiliary function (see, e.g., \cite{E} or  \cite{Pro}). 

First we set $M_1 := \Vert u\Vert_{L^{\infty}(\R\times (0,T))}$ 
and fix $x_0,\, t_0$, and $\delta$ that satisfy the inequalities $x_0 < 0$, $0<t_0<T$, and $\delta>0$.
Now we introduce an auxiliary function  in the form
$$
w(x,t):= u(x,t) - \sup_{x\in \R} a(x)
- \delta\left( \frac{t^{\alpha}}{\Gamma(1+\alpha)}
- \int^x_0 \frac{1}{q(\xi)} d\xi\right), \quad
x\in \R, \, 0<t<T
$$ 
and choose $L>0$ sufficiently large such that  the inequality
$$
\int^0_{-L} \frac{1}{q(\xi)} d\xi 
> \frac{M_1 + \sup_{x\in \R} a(x)}{\delta}, \quad -L<x_0<0  \eqno{(4.1)}
$$
holds true. 
Then, by (4.1), we get the following inequality 
$$
w(x,-L)= u(-L,t) - \sup_{x\in \R} a(x)
- \delta\left( \frac{t^{\alpha}}{\Gamma(1+\alpha)}
+ \int^0_{-L} \frac{1}{q(\xi)} d\xi\right)
$$
$$
\le M_1 + \sup_{x\in \R} \vert a(x) \vert
- \delta\int^0_{-L} \frac{1}{q(\xi)} d\xi < 0, \quad 
0<t<T.                 \eqno{(4.2)}
$$
Furthermore,
$$
w(x,0)= u(x,0) - \sup_{x\in \R} a(x)
+ \delta\int^x_0 \frac{1}{q(\xi)} d\xi
\le u(x,0) - \sup_{x\in \R} a(x) \le 0, \quad -L\le x\le 0,
$$
that is,
$$
w(x,0) \le 0, \quad -L\le x \le 0.   \eqno{(4.3)}
$$
On the other hand, direct calculations yield
$$
\pppa w(x,t) + q(x)\ppp_x w(x,t) = F(x,t) \le 0, \quad
-L < x < 0, \, 0<t<T.                                 \eqno{(4.4)}
$$
The inequalities (4.2)-(4.4) allow us to apply Theorem \ref{t1} (formula (2.4)) that leads  to the inequality
$$
w(x,t) \le 0, \quad -L\le x\le 0, \, 0\le t \le T
$$
and thus we arrive at the estimate 
$$
u(x_0,t) \le \sup_{x\in \R} a(x)
+ \delta\left( \frac{t^{\alpha}}{\Gamma(1+\alpha)}
- \int^{x_0}_0 \frac{1}{q(\xi)} d\xi\right), \quad 0<t<T.
$$ 
In the last formula, we let $\delta$ go to zero and get the inequality $u(x_0,t) \le \sup_{x\in \R} a(x)$.
Since the point $x_0<0$ is arbitrarily chosen, we have proved that
$$
u(x,t) \le \sup_{x\in \R} a(x) \quad \mbox{if } x\le 0.    \eqno{(4.5)}
$$
Introducing a new variable $y:= x+x_1$ with an arbitrarily chosen $x_1 > 0$, we can transfer the
previous arguments to any interval $(-L+x_1, x_1)$ and thus  arrive at the inequality
$u(x,t) \le \sup_{x\in \R} a(x)$ for $x\le x_1$ and $0<t<T$.
Since the point $x_1$ can be arbitrarily chosen, we have proved the relation (2.11).
Because $\sup v = -\inf (-v)$, the relation (2.12) can be derived from the relation (2.11)  by changing the signs in the equation (2.8) and in the initial condition (2.9) and  considering $-u$ instead of $u$.
The proof of Theorem \ref{t2} is completed.
\end{proof}
\section{Conclusions and directions for further research}
\label{s5}

%{\bf 5.1.}\\

In this paper, we proved a maximum principle for the general multi-term space-time-fractional transport equation and 
applied it for analysis of solutions to the 
initial-boundary-value problems for this equation. Here we restricted ourselves 
to the case of the one-dimensional fractional transport equation.  However, our arguments  can be transferred to the multi-dimensional case without any essential changes. Say, one  
can similarly treat the multi-term time-fractional transport equation
$$
\sump p_i(x,t)\ppp_t^{\alpha_i}u(x,t) 
- A(x,t)\cdot \nabla u(x,t) = r(x,t)u(x,t) + F(x,t), 
\quad x\in \OOO, \, 0<t<T,
$$
where $\OOO \subset \R^d$ is a bounded domain, $x = (x_1, ..., x_d) 
\in \R^d$, $A(x,t) = (a_1(x,t), ..., a_d(x,t))$, and
$\nabla v(x) = \left( \frac{\ppp v}{\ppp x_1}, ..., \frac{\ppp v}{\ppp x_d}
\right)$.  This equation will be considered elsewhere.

For validity of the results presented in this paper, we assumed that the zeroth order coefficient 
$r=r(x,t)$ of the fractional transport equation is non-positive on the whole domain $Q_T$. However, it is not clear if  this condition can be weakened or even removed. This problem is also a topic for our further research. 

In this paper, we did not address any nonlinear equations.  However, at least for some  semilinear 
fractional transport equations, our arguments still work and several important results can be derived.
For example, let us consider the equation
$$
\sump p_i(x,t)\ppp_t^{\alpha_i}u(x,t) 
+ \sumq q_j(x,t)\ppp_x^{\beta_j}u(x,t)
= r(x,t)u(x,t) + f(u(x,t))
                                                  \eqno{(5.1)}
$$
on the finite domain $Q_T:= (0,\ell) \times (0,T)$ and assume that the conditions (2.1) and (2.3) hold true.  Moreover, we suppose that the semilinear term  from the equation (5.1)  belongs to the following 
admissible set $\mathcal{F}$ of functions:
$$
\mathcal{F} := \left\{ f \in C^1(\R); \,
\frac{df}{d\xi}(\xi) \le 0,\ \xi \in \R\right\}.
$$
In fact,  the set $\mathcal{F}$ of admissible functions can be extended, but here we do not pursue 
the generality and prefer to focus on the underlying ideas.

In what follows, by $u_f = u_f(x,t)$ we denote a function that satisfies the inclusions (1.6) and the equation (5.1) with the semilinear term
$f$ from $\mathcal{F}$.  Then the following result holds true:
\\
{\it Let $f_1, f_2 \in \mathcal{F}$.
If 
$$
f_1(\xi) \le f_2(\xi), \quad \xi \in \R
$$
and
$$
u_{f_1}(0,t) \le u_{f_2}(0,t), \ \  0<t<T, \ \ 
u_{f_1}(x,0) \le u_{f_2}(x,0), \ \ 0<x<\ell, 
$$
then 
$$
u_{f_1}(x,t) \le u_{f_2}(x,t), \ \ (x,t) \in Q_T.
$$
}

Let us prove this statement. Setting $u:= u_{f_1} - u_{f_2}$, by the mean value theorem
we have the representation
$$
f_1(u_{f_1}(x,t)) - f_1(u_{f_2}(x,t)) 
= \frac{df_1}{d\xi}(\eta)u(x,t) =: g(x,t)u(x,t),
$$
where $\eta = \eta(x,t)$ is a number from the interval $[u_{f_1}(x,t),u_{f_2}(x,t)]$.  Because $f_1 \in \mathcal{F}$, the function $\eta= \eta(x,t)$ is a continuous function in both variables.  
Now we employ this representation and the identity
$$
f_1(u_{f_1}(x,t)) - f_2(u_{f_2}(x,t))
= (f_1(u_{f_1}(x,t)) - f_1(u_{f_2}(x,t)))
+ (f_1(u_{f_2}(x,t)) - f_2(u_{f_2}(x,t)))
$$
to rewrite the equation (5.1) as follows
\begin{align*}
& \sump p_i(x,t)\ppp_t^{\alpha_i}u(x,t) 
+ \sumq q_j(x,t)\ppp_x^{\beta_j}u(x,t)\\
=& (r(x,t)+g(x,t))u(x,t)
+ (f_1(u_{f_2}(x,t)) - f_2(u_{f_2}(x,t)), \
(x,t) \in Q_T.
\end{align*}
In the last equation,  $(r+g)(x,t) \le 0$, $(f_1(u_{f_2}(x,t)) - f_2(u_{f_2}(x,t))
\le 0$ for $(x,t) \in Q_T$,
$u(0,t)\le 0$ for $0\le t\le T$ and $u(x,0) \le 0$ for $0\le x \le \ell$. Thus we can apply Theorem \ref{t1} (the formula (2.4)) and obtain the inequality $u(x,t) \le 0$ for $(x,t) \in \qqqt$, that is,
$u_{f_1}(x,t) \le u_{f_2}(x,t)$ for $(x,t) \in Q_T$.

The last remark concerns Theorem \ref{t2} for the Cauchy problem (2.8)-(2.9) for the time-fractional transport equation. We state that the result formulated in Theorem \ref{t2} is  valid for a more general time-fractional transport equation in place of  the equation (2.8):
$$
\ppp_t^{\alpha_n}u(x,t) + \sum_{i=1}^{n-1} p_i(t)\ppp_t^{\alpha_i}u(x,t)
+ q(x)\ppp_x u(x,t) = F(x,t), \quad x\in \R, \, 0<t<T,
$$
where $0<\alpha_1 < ... < \alpha_n \le 1$ and $p_i(t) \ge 0$ for 
$0\le t \le T$ and $1 \le i \le n-1$.

In the rest of this section, we present a short sketch of its proof. The results presented in Chapter 3 of 
\cite{KRY} ensure existence and uniqueness of solution 
$u_0=u_0(t)$ to the initial-value problem
$$
\ppp_t^{\alpha_n}u_0(t) + \sum_{i=1}^{n-1} p_i(t)\ppp_t^{\alpha_i}u_0(t)
= 1, \quad 0<t<T, \quad u_0(0) = 0.
$$
Its solution $u_0 = u_0(t)$ is employed to define an auxiliary function in the form
$$
w(x,t) := u(x,t) - \sup_{x\in \R} a(x)
- \delta\left( u_0(t) - \int_0^x \frac{1}{q(\xi)} d\xi\right),
\quad x \in \R, \, 0<t<T.
$$
Now we suitably modify the condition (4.1), choose $L>0$ sufficiently large, and proceed as in the proof of Theorem \ref{t2} from  Section \ref{s4}. A complete version of the proof will be presented elsewhere.

\section*{Acknowledgment}
The third named author was supported by Grant-in-Aid for Scientific Research (S)
15H05740 of Japan Society for the Promotion of Science and
by the National Natural Science Foundation of China
(no. 11771270, 91730303).
This work was curried out with the support of the "RUDN University Program 5-100".

%%%%%%%%%%%%%%%%%%%%%%%%%%%%%%%%%%%%%%%%%%%%%%%%%
%%%%%%%%%%%%%%%%%%%%%%%%%%%%%%%%%%%%%%%%%%%%%%%%%

\end{document}